\documentclass[12pt]{article}

\def\@typesizes{%
       \or{5}{6.5}\or{6}{7.5}\or{7}{8.5}\or{8}{11}\or{9}{12}%
       \or{10}{13}
       \or{\@xipt}{14}\or{\@xiipt}{15}\or{\@xivpt}{18}%
       \or{\@xviipt}{20}\or{\@xxpt}{24}}

\usepackage{amsthm}
\usepackage{amsmath}
\usepackage{amssymb}
\usepackage{graphicx}
\usepackage{enumerate}
\allowdisplaybreaks
\numberwithin{equation}{section}
\numberwithin{figure}{section}

\theoremstyle{plain}
\newtheorem{theorem}{ Theorem}[section]
\newtheorem{proposition}[theorem]{ Proposition}
\newtheorem{lemma}[theorem]{ Lemma}
\newtheorem{corollary}[theorem]{ Corollary}
\newtheorem{example}[theorem]{ Example}
\newtheorem{remark}[theorem]{ Remark}
\newtheorem{definition}[theorem]{ Definition}
\newtheorem{conjecture}{ Conjecture}

\def\BET{\begin{theorem}}
\def\ENT{\end{theorem}}
\def\BEP{\begin{proposition}}
\def\ENP{\end{proposition}}
\def\BEL{\begin{lemma}}
\def\ENL{\end{lemma}}
\def\BEC{\begin{corollary}}
\def\ENC{\end{corollary}}
\def\BEE{\begin{example} \rm}
\def\ENE{\end{example}}
\def\BER{\begin{remark} \rm}
\def\ENR{\end{remark}}
\def\BED{\begin{definition} \rm}
\def\END{\end{definition}}
\def\BECJ{\begin{conjecture}}
\def\ENCJ{\end{conjecture}}

\def\bea{\begin{eqnarray}}
\def\eea{\end{eqnarray}}

\def\beq{\begin{equation}}
\def\eeq{\end{equation}}

\def\beal{\begin{align*}}

\def\eeal{ \end{align*} }

\begin{document}

\begin{center}
{\large{\bf Criteria for the Absence and Existence of Bounded Solutions at the Threshold Frequency in a Junction of Quantum Waveguides}}
\end{center}

\bigskip

\centerline{F. L. Bakharev$^{a,b}$ and S. A. Nazarov$^{a,c}$}

\bigskip

\begin{center}

\emph{a) St. Petersburg State University, Mathematics and Mechanics Faculty, 7/9 Universitetskaya nab., St. Petersburg, 199034 Russia} \\
\emph{b) St. Petersburg State University, Chebyshev Laboratory, 14th Line V.O., 29B, Saint Petersburg 199178 Russia} \\
\emph{c) Institute of Problems of Mechanical Engineering RAS, V.O., Bolshoj pr., 61, St. Petersburg, 199178 Russia} \\

\texttt{fbakharev@yandex.ru, f.bakharev@spbu.ru, srgnazarov@yahoo.co.uk}

\end{center}

{\small

\noindent {\bf Abstract}: 
In the junction $\Omega$ of several semi-infinite cylindrical waveguides we consider the Dirichlet
Laplacian whose continuous spectrum is the ray $[\lambda_\dagger, +\infty)$ with a positive cut-off
value $\lambda_\dagger$. We give two different criteria for the threshold resonance 
generated by nontrivial bounded solutions to the Dirichlet problem
for the Helmholtz equation $-\Delta u=\lambda_\dagger u$ in $\Omega$. The first criterion
is quite simple and is convenient to disprove the existence of bounded solutions. The
second criterion is rather involved but can help to detect concrete
shapes supporting the resonance. Moreover, the latter distinguishes in a natural way
between stabilizing, i.e., bounded but non-descending solutions and trapped modes with
exponential decay at infinity. 

\medskip


}

\section{Introduction}

\subsection{Motivation}
In a domain with several cylindrical outlets to infinity, Fig. 1, we are interested in retrieving the threshold
resonance generated by nontrivial bounded solutions of the spectral Dirichlet problem for the Laplace operator when the
spectral parameter coincides with the lower bound $\lambda_\dagger$ of the continuous spectrum. This concern is caused by the dimension
reduction procedure for lattices of thin waveguides, namely, according to \cite{Gri, MoVa}, transmission conditions
at the vertices of the graph skeleton in the one dimensional model of the lattice crucially depend on whether the boundary-value problem in the stretched node, 
Fig.2, admits stabilizing (bounded but not decaying) solutions to the homogeneous Dirichlet problem. 
For acoustic waveguides with hard walls, cf. \cite{KuchChen, ExPo}, the Neumann problem for the Laplace operator  
surely gets such solutions, namely constants (the threshold is null). For quantum waveguides described by the Dirichlet problem, 
the existence and absence questions are much more delicate because of the positive threshold $\lambda_\dagger>0$. 
Certain sufficient conditions \cite{Gri, Pank} and concrete
canonical shapes \cite{na566, na578, na602, naIZV} are known to assure the absence of bounded solutions at the threshold. 
At the same time, as was indirectly verified in \cite{na566, na602, naIZV}, bounded solutions may 
emerge in parameter dependent
junctions but only at isolated values of the inserted geometrical parameter. 

\begin{figure}[t]
\begin{center}
\includegraphics[scale=0.4]{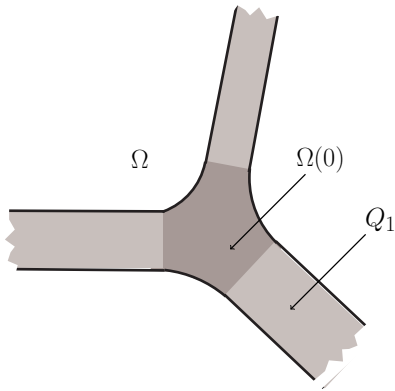}
\end{center}
\caption{The junction of semi-infinite quantum waveguides}
\end{figure}

In this paper we present two quite different criteria for the threshold resonance and distinguish between them with
the following reason. The first criterion in Section 2 with rather simple formulation is convenient to verify
the absence of bounded solutions at the lower bound of the continuous spectrum but we do not see a way to apply 
this criterion to finding a particular bounded solution in a specific geometry. On the contrary, the second criterion in
the Section 3 requiring for several definitions of auxiliary objects, can be employed to develop
analytical, in particular, asymptotic methods or numerical schemes to 
detect and analyse concrete stabilizing (i.e. bounded but non-decaying) solutions and trapped modes
with the exponential decay at infinity. At the same time, these methods and schemes may also help to disprove the
threshold resonance but the latter is much more expensive in comparison with the first, absence, criterion.

Our proofs in Section 2 are conducted in such a way that they can be easily adapted for other problems, e.g., for 
mixed boundary conditions \cite{Vas}. Nevertheless, any generalization
of the whole existence criterion in Section 3 is still a fully open question.

\begin{figure}[t]
\begin{center}
\includegraphics[scale=0.5]{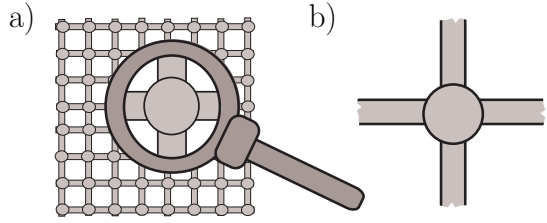}
\end{center}
\caption{Rectangular lattice (a) and the infinite cruciform waveguide (b) to describe the boundary layer phenomenon}
\end{figure}

\subsection{Statement of the spectral problem}
Let $\Omega\subset \mathbb{R}^d$, $d\geq 2$ be a domain, by definition an open connected set, with several cylindrical
outlets $Q_1$, \ldots, $Q_N$ to infinity. We assume that 
$\Omega=\Omega(0)\cup Q_1\cup\ldots\cup Q_N$ where $\Omega(0)$ 
is a bounded domain (shaded in Fig. 1) with Lipschitz boundary, and 
$Q_n\cap Q_k=\varnothing$ for $n\ne k$, $Q_n\cap \Omega(0)=\varnothing$ for $n=1,\ldots, N$.
In each outlet $Q_n=\omega_n\times [0,\infty)$ we introduce the Cartesian system $(y_n,z_n)$ 
of local coordinates  with $y_n\in\omega_n$, $z_n\in[0,+\infty)$, where the cross-section 
$\omega_n\subset \mathbb{R}^{d-1}$ is a bounded domain with Lipschitz boundary 
$\partial \omega_n$. Note that the outlets $Q_n$ include their ends, i.e. $\omega_n\times \{0\}\subset \partial\Omega(0)$ for $n=1,\ldots,N$.  
We also will deal with the truncated waveguide   
\begin{equation}
\label{1.0}
\Omega(R)=\Omega(0) \cup \bigcup_{n=1}^N Q_n(R), \quad Q_n(R)=\{x\in Q_n : z_n\in[0,R)\}, \quad R>0.
\end{equation}
In what follows we use the notation
$$
Q_n^R=Q_n \setminus \overline{Q_n(R)}=\{x\in Q_n\colon z_n>R\}, \quad R>0, 
$$
and the index $n$ is usually omitted in proofs related to any outlet.

We consider the spectral problem for the Laplacian $\Delta=\nabla\cdot\nabla$
\begin{equation}
\label{1.1}
-\Delta u =\lambda u \mbox{ in } \Omega, \quad u=0 \mbox{ on } \partial \Omega, 
\end{equation}  
where $\nabla=\mathop{\rm grad}\nolimits$, $\lambda$ is a spectral parameter and $\partial \Omega$ is the boundary of 
$\Omega$ which, for simplicity, is assumed to be Lipschitz. 

The variational form of the problem \eqref{1.1} reads:
\begin{equation}
\label{1.2}
(\nabla u,\nabla \psi)_{\Omega}=\lambda (u,\psi)_{\Omega} \quad \forall \psi\in H^1_0(\Omega),
\end{equation} 
where $(\cdot,\cdot)_\Omega$ is the natural scalar product in the Lebesgue space $L^2(\Omega)$
and $H^1_0(\Omega)$ stands for the Sobolev space of functions
vanishing at the boundary $\partial \Omega$. Since the bilinear form on the left-hand side of the integral identity \eqref{1.2} is
closed and positive definite, it gives rise \cite[Ch.10]{BiSo}, \cite[Ch.VIII]{RS1} to an unbounded positive definite
self-adjoint operator $\mathcal{A}$ in the Hilbert space $L^2(\Omega)$.

The Dirichlet problem on the cross-section
\begin{equation}
\label{A2}
-\Delta_y \Phi^n(y)=\Lambda^n\Phi^n(y), \, y\in\omega_n,\quad \Phi^n(y)=0,\, y\in\partial\omega_n,
\end{equation}
has the monotone and unbounded eigenvalue sequence
$$
0<\Lambda_1^n<\Lambda_2^n\leq\Lambda_3^n\leq\ldots\leq \Lambda_j^n \leq \ldots \to+\infty
$$
and the corresponding real eigenfunctions $\Phi_j^n\in H^1_0(\omega_n)$, $j\in\mathbb{N}$ are subject to the orthogonality 
and normalization conditions 
\begin{equation}
\label{A4}
(\Phi_j^n,\Phi_k^n)_{\omega_n}=\delta_{j,k},\quad j,k\in\mathbb{N},
\end{equation}
where $\delta_{j,k}$ is the Kronecker symbol and ${\mathbb N}=\{1,2,3,\ldots\}$.

It is known that the continuous spectrum $\sigma_c$ of the operator $\mathcal{A}$ is the ray 
$[\lambda_\dagger,+\infty)$ where the lower bound $\lambda_\dagger=\min \{\Lambda_1^1,\ldots,\Lambda _1^N\} >0$
coincides with the smallest among all principal (minimal) eigenvalues $\Lambda_1^n$.
The total multiplicity $\varkappa:=\# \sigma_d$ of the discrete spectrum 
\begin{equation}
\label{1.3}
0<\lambda_1<\lambda_2\leq\ldots\leq \lambda_{\varkappa}<\lambda_\dagger
\end{equation}
of the operator $\mathcal{A}$ is known to be finite.

If a cranked waveguide belongs to $\Omega$ and is composed of two skewed
semi-infinite cylinders which have the cross-sections congruent to $\omega_1$  
and meet each other under the angle $\alpha\in (0,\pi)$,
then $\varkappa\geq 1$ according to a result in \cite{AvBeGiMa} and the max-min principle \cite[Th 10.2.2]{BiSo}. Furthermore,
the papers \cite{AvBeGiMa, na561, Dauge, Dauge1} and \cite{BaMaNaZaa} give examples of arbitrary large $\varkappa$ in
dimension 2 and 3, respectively. We refer the book \cite{ExnerBook} for a completed review of results on the discrete
spectrum of quantum waveguides and their junctions.

\subsection{Trapped modes and stabilizing solutions}
Within the approach \cite{Gri, MoVa}, it is important to distinguish between stabilizing 
solutions and trapped modes. To explain the main difference between these kinds of 
bounded solutions, we consider a thin, of diameter $\varepsilon\ll 1$,
finite lattice $\Upsilon^\varepsilon$ of quantum waveguides and its fragment
$$
\Upsilon^\varepsilon_\bullet=\{x\colon \varepsilon^{-1}(x-x_\bullet)\in \Omega(1/\varepsilon)\}
$$
around the node $\upsilon^\varepsilon_\bullet=\{x\colon \varepsilon^{-1}(x-x_\bullet)\in\Omega(0)\}$ with the center $x_\bullet$.
To simplify formulas, we suppose for a while that all cylinders $Q_1, \ldots, Q_N$ have the same cross-section $\omega$
of unit $(d-1)$-dimensional area. If in addition to the isolated eigenvalues \eqref{1.3}, the operator $\mathcal{A}$
has the embedded eigenvalue $\lambda_\dagger$ of multiplicity ${{\bf k}}\geq 0$, then, according
to \cite{Gri}, the Dirichlet problem in $\Upsilon^\varepsilon$ gets eigenvalues with the asymptotic forms 
\begin{equation}
\label{M1}
\begin{array}{l}
M_j^\varepsilon=\varepsilon^{-2}\lambda_j+o(e^{-\delta_j/\varepsilon}),\quad \delta_j>0, \quad j=1,\ldots, \varkappa,\\
M_{k+\varkappa}^\varepsilon=\varepsilon^{-2}\lambda_\dagger + o(e^{-\delta_\dagger/\varepsilon}), \quad \delta_\dagger>0,\quad  k=1,\ldots,{\bf k}.
\end{array}
\end{equation} 
The corresponding eigenfunctions are localized in the vicinity of the node $\upsilon^\varepsilon_\bullet$
and become exponentially small at a distance from it.

Stabilizing solutions in $\Omega$ at the threshold $\lambda=\lambda_\dagger$ influence the spectrum in $\Upsilon_\bullet^\varepsilon$
in a quite different way. Indeed, eigenvalues above the rescaled, cf. \eqref{M1}, threshold $\varepsilon^{-2}\lambda_\dagger$
are determined through ordinary differential equations on edges of the skeleton $\Upsilon^0=\cap_{\varepsilon>0}\Upsilon^\varepsilon$
linked by certain transmission conditions at vertices of the graph $\Upsilon^0$. If the problem in the infinite waveguide \eqref{1.0}
has no stabilizing solutions at the threshold, then the transmission conditions at the vertex $x_\bullet$ are nothing but the Dirichlet
ones, i.e. eigenfunctions in the one-dimensional model must vanish at this vertex and, therefore, the graph
edges emerging from $x_\bullet$ decouple.
On the other hand, according to \cite{Gri, MoVa}, the existence of stabilizing solutions changes the Dirichlet conditions at $x_\bullet$
for some other conditions, in particular, the Kirchhoff ones like in the Pauling model \cite{Pau} for the
Neumann problem \cite{KuchChen, ExPo}. Thereby, the main question in the framework of the dimension reduction procedure 
\cite{Gri, MoVa} becomes to detect
stabilizing solutions rather than all bounded solutions and the corresponding threshold resonance. The existence criterion in Section 3
makes the necessary separation of two kinds of bounded solutions in a natural way, compare Proposition \ref{Pr1} and Proposition 3.1. However, 
the absence criterion in Section 2 cannot directly select stabilizing solution and we provide in Section 2.4 a simple sufficient condition
for absence of trapped modes but do not know an appropriate necessary condition yet.  

\section{An absence criterion}

\subsection{Formulation of the first criterion}
We consider the auxiliary spectral problem with mixed  boundary conditions
\begin{equation}
\label{1.4}
\begin{array}{l}
-\Delta v^R=\mu^R v^R \mbox{ in } \Omega(R),\quad v^R=0 \mbox{ on } \Gamma(R):=\partial\Omega(R)\cap \partial\Omega,
\\
\partial_\nu v^R=0 \mbox{ on } \gamma(R)=\partial\Omega(R) \setminus\partial\Omega,
\end{array}
\end{equation} 
where $R\geq 0$, $\partial_\nu$ is the outward normal derivative, in particular, $\partial_\nu=\partial_{z_n}=\partial/\partial z_n$ 
on the truncation surface $\gamma_{n}(R)=\{x \colon y_n\in\omega_n, z_n=R\}$.

The variational formulation of the problem \eqref{1.4} reads:
\begin{equation}
\label{1.6}
(\nabla v^R, \nabla \psi)_{\Omega(R)} =\mu^R (v^R,\psi)_{\Omega(R)} \quad \forall \psi\in {\mathcal H}^0_R
\end{equation}
where ${\mathcal H}^0_R:=H^1_0(\Omega(R),\Gamma(R))$ is a subspace of functions in $H^1(\Omega(R))$ vanishing at $\Gamma(R)$.
The problem \eqref{1.6}  gives rise to unbounded positive definite and self-adjoint operator $\mathcal{A}^R$ in 
$L^2(\Omega(R))$. Since ${\mathcal H}^0_R$ is compactly embedded into $L^2(\Omega(R))$, the spectrum of $\mathcal{A}^R$
is discrete and composes the monotone unbounded sequence of eigenvalues
\begin{equation}
\label{mu-mu}
0<\mu_1^R<\mu_2^R\leq \ldots \leq \mu_\varkappa^R \leq \mu_{\varkappa+1}^R\leq \ldots \to+\infty\,,
\end{equation}
where their multiplicity is taken into account. We will prove the following criterion for the threshold resonance. 

\begin{theorem}
\label{Crit1}
The problem \eqref{1.1} has no threshold resonance if and only if, for some $R\geq 0$, the eigenvalue 
$\mu_{\varkappa+1}^R$ of the problem \eqref{1.4} meets the inequality  $\mu_{\varkappa+1}^R>\lambda_\dagger$. 
Here, $\varkappa$ 
is the total multiplicity of the discrete spectrum $\sigma_d$, see~\eqref{1.3}.
\end{theorem}

\subsection{Sufficiency}

\begin{proposition}
\label{Pr2.2}
If the eigenvalue $\mu_{\varkappa+1}^R$ of the problem \eqref{1.4} meets the inequality $\mu_{\varkappa+1}^R>\lambda_\dagger$
for some $R\geq 0$, then the threshold resonance is absent in the problem \eqref{1.1}.
\end{proposition}

This result coincides with Theorem 3 in \cite{Pank}. Here, we only provide a short sketch of a proof. 
The proof is based on a
simple observation, originally used in \cite{na566, na578, na602, naIZV} for the case $\varkappa=1$: 
if the threshold resonance occurs, one may construct a small compact perturbation $\mathcal{B}$
located in $\Omega(0)$, that is $\mathcal{B} u=0$ in $\Omega\setminus \Omega(0)$, with the following
properties. First of all, the perturbed eigenvalues  $\widehat\mu_1^R,\ldots, \widehat\mu_\varkappa^R$
and $\widehat\mu_{\varkappa+1}^R$ of the operator $\widehat{\mathcal{A}}_R=\mathcal{A}_R+\mathcal{B}$ still stay, respectively, below
and above the threshold $\lambda_\dagger$, so that one gets the Poincare inequality
\begin{equation}
\label{poin}
(\widehat{\mathcal{A}}_R u_0,u_0)_{\Omega(R)}\geq \widehat{\mu}^R_{\varkappa+1}\|u_0;L^2(\Omega(R))\|^2\geq \lambda_\dagger \|u_0;L^2(\Omega(R))\|^2,
\end{equation}
where $u_0$ is orthogonal in $L^2(\Omega(R))$ to eigenfunctions of $\widehat{\mathcal{A}}_R$ corresponding to
$\widehat\mu_{1}^R, \ldots, \widehat\mu_{\varkappa}^R$.
Then in a standard way the max-min principle, cf. \cite[Theorem 10.2.2]{BiSo}, equipped with the inequalities \eqref{poin}
and
\begin{equation}
\label{fri}
\|\nabla u; L^2(Q_n^R)\|^2\geq \Lambda_1^n \|u; L^2(Q_n^R)\|^2\geq\lambda_\dagger\|u;L^2(Q_n^R)\|^2, 
\end{equation}
verifies that the total multiplicity of the discrete spectrum $\widehat{\sigma}_d$
of the operator $\widehat{\mathcal{A}}=\mathcal{A}+\mathcal{B}$ meats the inequality $\#\widehat{\sigma}_d\leq \varkappa$. 
Notice that \eqref{fri} is a direct consequence of the Friedrichs inequality in the cross-section~$\omega_n$. 

Finally, a special choice, see cf. \cite{na566, na578, na602, naIZV, Pank},  
of the perturbation $\mathcal{B}$ provides the existence of the eigenvalue 
$\widehat{\mu}_{\varkappa+1}<\lambda_\dagger$ in the 
discrete spectrum $\widehat{\sigma}_d$ of the perturbed operator~$\widehat{\mathcal{A}}$. 
The latter contradiction completes the proof of Proposition~\ref{Pr2.2}.

\subsection{Necessity}
We proceed with proving that eigenvalues in the sequence \eqref{mu-mu} below the continuous spectrum are monotone 
increasing functions in $R$.

\begin{lemma}
If the eigenvalue $\mu_{k}^R$ of the problem \eqref{1.4} meets the inequality $\mu_{k}^R<\lambda_\dagger$
for some $R>0$,
then there exists $r_0>0$ such that 
$$
\mu_{k}^R<\mu_{k}^{R+r}<\lambda_\dagger \quad \forall r\in(0,r_0)\,.
$$ 
\end{lemma}

\noindent{\bf Proof.}
We consider the operator $\mathcal{A}^{R+r}$ in $L^2(\Omega(R+r))$ for small $r>0$ as a perturbation of  $\mathcal{A}^{R}$
in a certain sense.
For the simple eigenvalue $\mu^R$ (we omit the index $k$), we denote by $v^R$ the corresponding 
eigenfunction normalized in $L^2(\Omega(R))$. 
Let us accept the simplest asymptotic ans\"atze 
\begin{equation}
\label{5.1}
\mu^{R+r}=\mu^R+r\mu'+\ldots,
\end{equation} 
\begin{equation}
\label{5.2}
v^{R+r}=v^R+r v'+\ldots
\end{equation}
where the correction terms $\mu'$ and $v'$ are to be determined and ellipses replace small reminders to 
be estimated.
The functions $v^R$ and $v'$ defined in $\Omega(R)$, can be smoothly extended onto $\Omega \supset \Omega(R+r)$.
We use the same letters for these extensions. Plugging formulas \eqref{5.1} and \eqref{5.2} into 
the equation for $v^{R+r}$ on $\Omega(R)$ and collecting terms
of the same order in  $r$ yield
\begin{equation}
\label{5.3}
\Delta v'(x)+\mu^R v'(x)=-\mu' v^R(x), \, x\in \Omega(R)\,.
\end{equation} 
Imposing the Dirichlet condition 
$$
v'(x)=0, \, x\in  \Gamma(R),
$$
is quite evident. The Neumann condition on 
$\gamma(R+r)$ can be formally transferred to $\gamma(R)$ by the Taylor formula in the variable $z$, indeed, 
\begin{multline*}
\partial_z v^{R+r}\big|_{z=R+r}=\partial_z v^R \big|_{z=R+r}+r \partial_z v'\big|_{z=R+r}+\ldots=\\
=
\partial_z v^R\big|_{z=R}+r\partial^2_z v^R \big|_{z=R}+r \partial_z v'\big|_{z=R}+\ldots
\end{multline*}
We recall the Helmholtz equation for $v^R$ and introduce the boundary condition 
\begin{equation}
\label{5.4}
\partial_z v'(x)=-\partial_z^2 v^R(x)=\Delta_y v^R(x) +\mu^R v^R(x), \, x\in \gamma(R)\,.
\end{equation}

The compatibility condition in the problem \eqref{5.3}-\eqref{5.4} reads:
\begin{multline*}
\mu'=\mu'(v^R,v^R)_{\Omega(R)}=-({\Delta}_y v^R +\mu^R v^R, v^R)_{\gamma(R)}=\\
=\|{\nabla}_y v^R;L^2(\gamma(R))\|^2-\mu^R \|v^R; L^2(\gamma(R))\|^2\,.
\end{multline*}
By the Friedrichs inequality on $\gamma(R)$, we obtain 
$$
\mu'\geq \sum_{n=1}^N (\Lambda_1^n-\mu^R)\|v^R; L^2(\gamma_n(R))\|^2\geq (\lambda_\dagger-\mu^R)\|v^R; L^2(\gamma(R))\|^2>0.
$$ 

If the eigenvalue $\mu^R$ has multiplicity $m$, calculations mainly remain the same. The leading term in the anzatz \eqref{5.2}
becomes a linear combination of the corresponding eigenfunctions $v^{R}_1$, $v^{R}_2$, \ldots, $v^{R}_m$ 
orthonormalized in $L^2(\Omega(R))$ with the coefficient column ${\bf c}=(c_1,c_2,\ldots, c_m)^\top$. 
Repeating the above calculations  with minor modifications, 
we observe that the correction terms $\Lambda_1'$, \ldots, $\Lambda_m'$ in \eqref{5.1}
are found from the system of linear algebraic equations 
$$
{\bf A} {\bf c}=\Lambda' {\bf c}, 
$$
where the self-adjoint and positive definite matrix ${\bf A}$ of size  $m\times m$ has the entries 
$$
{\bf A}_{pq}=(\nabla_y v^{R}_p, \nabla_y v^{R}_q)_{\gamma(R)}-\mu^R (v^{R}_p, v^{R}_q)_{\gamma(R)}, \, p,q=1,\ldots,m.
$$

The correction terms in the asymptotic formula \eqref{5.1} for the eigenvalues $\mu_k^{R+r},\ldots, \mu_{k+m-1}^{R+r}$
of the problem \eqref{1.4} in $\Omega(R+r)$ involve the eigenvalues $\Lambda_1'$,\ldots, $\Lambda_m'$
of ${\bf A}$ and therefore become strictly positive as in the case of simple eigenvalues. To conclude
with the proof, we mention that the error estimates 
$|\mu_k^{R+r}-\mu_k^R-r\mu_k'|\leq c_k r^2$ for $r\in (0,r_k)$ 
are derived in a classical way, see \cite[Ch.7, \S6.5]{Kato},
because one can readily construct ``almost identical'' diffeomorphism between the domains $\Omega(R+r)$ and $\Omega(R)$, 
which is identical inside $\Omega(R-1)$ and coincides with the shift operator near the faces $\gamma_n(R)$.
We omit here the corresponding simple and routine computations.
\hfill$\Box$

\medskip

Now assume that the condition on $\mu_{\varkappa+1}^R$ in Theorem 2.1 is violated.
This means that, in particular, $\mu^R=\mu_{\varkappa+1}^R < \lambda_\dagger$ for all $R>2$.
We normalize the corresponding eigenfunction $v_R$
as follows: 
\begin{equation}
\label{normal}
\| v^R;L^2(\Omega(2))\|=1. 
\end{equation}

We are going to verify that there exists a monotone unbounded sequence $\{R_j\}_{j\in \mathbb{N}}$ such that $v^{R_j}$
converges in a certain sense to a non-trivial bounded solution $v^\infty$ of the problem \eqref{1.1}
with parameter $\mu^\infty=\lim \mu^{R_j}$ as $j\to+\infty$. To this end, we use the decomposition 
\begin{equation}
\label{Lm2.4-1}
\chi_1(x)v^R(x)=\chi_1(x)w_{n}^R(z_n) \Phi_1^n(y_n)+v^R_\perp(x) \quad \mbox{in}\quad Q_n, \quad n=1,\ldots, N
\end{equation}
and treat its ingredients $w_n^R$ and $v_{\perp}^R$ in a different way.

Let us recall that $\Phi_1^n$ is the first eigenfunction of the Dirichlet Laplacian in $\omega_n$ 
and $\|\Phi_1^n;L^2(\omega_n)\|=1$. Furthermore,
\begin{equation}
\label{ort1}
w_n^R(z_n)=\int_{\omega_n} v^R(x) \Phi_1^n(y_n) dy_n,\quad
\int_{\omega_n} v^R_\perp(y_n,z_n) \Phi_1^n(y_n)dy_n=0.
\end{equation} 
The smooth cut-off function $\chi_s$ is chosen such that $0\leq \chi_s \leq 1$ and
\begin{eqnarray*}
&&\chi_s(x)=0 \mbox{ if } x\in\Omega(s-1),
\quad \chi_s(x)=\chi(z_n-s)  \mbox{ if }  z_n\in [s-1,s], \\
&&\chi_s(x)=1  \mbox{ if }  z_n\geq s,
\end{eqnarray*}
with a fixed smooth  function $\chi$. We will also use the difference $\mathcal{X}_s(x)=1-\chi_s(x)$.
We further define $w^R$ as follows: 
$$\chi_1 w^R=\chi_1 w_n^R\Phi_1^n  \mbox{ in } Q_n,\quad w^R=0 \mbox{ in } \Omega(0).$$
Note that that $v^R_\perp$ in \eqref{ort1} is assumed to be zero in $\Omega(0)$.

\begin{lemma}
There exists a positive constant $c_1(\Omega)$ such that
\begin{equation}
\label{l1}
\|\nabla_x v^R;L^2(\Omega(1))\|+\sum_{n=1}^{N}\|v^R;L^2(\gamma_n(1))\|\leq c_1(\Omega)\,.
\end{equation}
\end{lemma}

\noindent{\bf Proof.}
From the integral identity \eqref{1.6} we derive the relation
\begin{multline*}
\|\nabla v^R;L^2(\Omega(1))\|^2\leq \|\nabla (\mathcal{X}_2 v^R);L^2(\Omega(2))\|^2=\\=
\mu^R\|\mathcal{X}_2 v^R;L^2(\Omega(2))\|^2+\|v^R\nabla \mathcal{X}_2;L^2(\Omega(2))\|^2\leq c_1(\Omega)\,.
\end{multline*}
The last inequality follows from \eqref{normal}. The standard trace inequality provides the desired estimate 
of the norm $\|v^R;L^2(\gamma_n(1))\|$ as well. 
\hfill$\Box$

\bigskip

Separation of variables
gives
\begin{equation}
\label{Lm2.4-2}
-\partial^2_{z_n} w_n^R(z_n)=(\mu^R-\Lambda_1^n)w_n^R(z_n)  \mbox{ for }   z_n>1,\quad
\partial_{z_n}w_n^R(R)=0\,.
\end{equation}
Moreover, formulas \eqref{Lm2.4-1} and \eqref{l1} assure that 
$|w_n^R(1)|\leq c_1(\Omega)$. 
A solution of the problem \eqref{Lm2.4-2}
takes the form
$$
w_n^R(z_n)=a_n^R (e^{-\alpha_n(R) z_j}+e^{-2\alpha_n(R) R}e^{\alpha_n(R)z_n})
$$
where $\alpha_n(R)=\sqrt{\Lambda_1^n-\mu^R}$. Thus,  
\begin{equation}
\label{aj}
|a_n^R|\leq c_2(\Omega), \quad n=1,\ldots, N, \,\, R>2.
\end{equation} 

Now we examine the function $v_\perp^R$ in \eqref{Lm2.4-1}. 
First, the Poincare inequality 
\begin{equation}
\label{2.4-4}
\|\nabla_{y_n} v^R_\perp(\cdot, z_n); L^2(\omega_n)\|^2\geq \Lambda_2^n\|v^R_\perp(\cdot,z_n);L^2(\omega_n)\|^2, \quad z_n>0,
\end{equation}
is valid due to the orthogonality condition in \eqref{ort1}.
Furthermore, $v^R_\perp$ is a solution of the problem
\begin{eqnarray}
\label{2.4-5}
&&-\Delta v^R_\perp -\mu^R v^R_\perp=[\Delta,\chi_1](v^R-w^R)=:f^R   \mbox{ in }  \Omega(R),\\
\label{2.4-6}
&&v^R_\perp=0   \mbox{ on }   \Gamma(R),\quad \partial_n v^R_\perp =0   \mbox{ on }   \gamma(R),
\end{eqnarray}
where $[\Delta,\chi_1]$ is the commutator of the Laplacian and the cut-off function $\chi_1$ (a first-order differential operator).
Obviously,
\begin{equation}
\label{26}
\mathop{\rm supp}\nolimits f^R \subset \Omega(1)\setminus \Omega(0),\quad \|f^R; L^2(\Omega(R))\| \leq c_3(\Omega).
\end{equation}

We fix a parameter $\beta=\beta(\Omega)$ such that 
\begin{equation}
\label{beta}
0<\beta <\frac{1}{2}\left(\min_{1\leq n\leq N}\Lambda_2^n-\lambda_\dagger\right)
\end{equation}
and introduce the weight function $T_\beta$,
$$
T_\beta(x)=1 \mbox{ for } x\in \Omega(1), \quad T_\beta(x)=e^{\beta(z_n-1)} \mbox{ for } z_n\geq 1\,.
$$
We also need the weighted Sobolev and Lebesgue spaces $W^1_\beta(\Omega)$ and $L^2_\beta(\Omega)$ with 
the following norms:
$$
\|v;W^1_\beta(\Omega)\|=\|T_\beta v; H^1(\Omega)\|  \quad \mbox{and} \quad \|v;L^2_\beta(\Omega)\|=\|T_\beta v; L^2(\Omega)\|. 
$$
If $\Omega$ is replaced with $\Omega(R)$ in these definitions, we obtain
the spaces $W^1_\beta(\Omega(R))$ and $L^2_\beta(\Omega(R))$ which coincide algebraically and topologically 
with $H^1(\Omega(R))$ and $L^2(\Omega(R))$,
respectively.

\begin{lemma}
For all $R>1$, the function $v_\perp^R$ enjoys the estimate
\begin{equation}
\label{vRbot}
\|v_\perp^R;W^1_\beta(\Omega(R))\|\leq c_\perp(\Omega).
\end{equation}
\end{lemma}

\noindent{\bf Proof.}
The function $T_\beta^2v_\perp^R$ falls into the space $\mathcal{H}_R^0$ and can be inserted as a test function into the 
integral identity for the problem \eqref{2.4-5}-\eqref{2.4-6}. Thus, we have
\begin{equation}
\label{2.21}
(\nabla v^R_\perp, \nabla(T_\beta^2 v^R_\perp))_{\Omega(R)}=\mu^R \|T_\beta v^R_\perp; L^2(\Omega(R))\|^2+(f^R,T^2_\beta v^R_\perp)_{\Omega(R)}\,.
\end{equation}
The left-hand side is equal to 
\begin{equation}
\label{2.22}
(\nabla v^R_\perp, \nabla(T_\beta^2 v^R_\perp))_{\Omega(R)}=
\|\nabla (T_\beta v^R_\perp);L^2(\Omega(R))\|^2-\|v^R_\perp \nabla(T_\beta);L^2(\Omega(R))\|^2
\end{equation}
and, in view of \eqref{2.4-4} and \eqref{beta}, gets the below bound
$$
\left(\min_{1\leq n\leq N}\Lambda_2^n-\beta\right) \|T_\beta v^R_\perp ;L^2(\Omega(R))\|^2\geq 
(\beta+\lambda_\dagger)\|T_\beta v^R_\perp ;L^2(\Omega(R))\|^2 .
$$
Hence, we deduce that
\begin{multline}
\label{2.23}
\|T_\beta v^R_\perp ;L^2(\Omega(R))\|^2 \leq \beta^{-1} (T^2_\beta f^R, v^R_\perp)_{\Omega(R)}\leq\\
\leq 
\beta^{-1} 
\|T^2_\beta f^R;L^2(\Omega(R))\| \|v^R_\perp; L^2(\Omega(1))\|\,.
\end{multline}
Relations \eqref{26}, \eqref{2.21}-\eqref{2.23} show that the product 
$T_\beta v^R_\perp$ also
enjoys the inequality
$$
\|\nabla (T_\beta v^R_\perp) ;L^2(\Omega(R))\|^2 \leq c_4(\Omega)
$$
with some constant $c_4(\Omega)$ and, therefore, the inequality \eqref{vRbot} holds true.
\hfill$\Box$

\medskip

Now, for $x\in \Omega(R)$, we determine the function 
\begin{equation}
\label{2.25}
\widehat{v}^R=\mathcal{X}_R \big(v^R-\chi_1w^R\big)=
\mathcal{X}_1v^R+\mathcal{X}_Rv^R_\perp\,, 
\end{equation}
and extend it by zero onto the whole domain $\Omega$. First of all,
\begin{equation}
\label{2.26}
\|\widehat{v}^R; W^1_\beta(\Omega)\|\leq \widehat{c}(\Omega).
\end{equation} 
The equation
\begin{equation}
\label{2.27}
-\Delta \widehat{v}^R-\mu^R \widehat{v}^R= [\Delta, \chi_1] w^R+[\Delta, \mathcal{X}_R]v^\perp=:g^R+h^R
\end{equation}
in the variational form becomes
\begin{equation}
\label{2.28}
( \nabla\widehat{v}^R, \nabla \psi)_\Omega - \mu^R (\widehat{v}^R, \psi)_\Omega= (g^R,\psi)_{\Omega}+(h^R,\psi)_{\Omega} \quad 
\forall \psi\in C^\infty_0(\Omega).
\end{equation}

We are going to perform the limit passage $R\to +\infty$ in \eqref{2.28}. Since $\mu^R$ is non-decreasing 
function in $R$, it has a limit, 
$$\lim_{R\to \infty}\mu^R=\mu^\infty\leq \lambda_\dagger.$$
The relations \eqref{2.26} and \eqref{aj} allows us to find a monotone unbounded sequence $\{R_k\}$ such that
\begin{eqnarray}
&& \widehat{v}^{R_k} \to \widehat{v}^{\,\infty} \quad \mbox{weakly in} \quad W^1_\beta(\Omega), \nonumber \\
\label{conv2}
&& \widehat{v}^{R_k} \to \widehat{v}^{\,\infty} \quad  \mbox{strongly in} \quad L^2(\Omega(2)),\\
\label{conv3}
&& a_n^{R_k} \to a_n^\infty, \,\, \alpha_n(R)\to \alpha_n^\infty,\,\, e^{-2\alpha_n(R)R}\to c_n^\infty   \mbox{ for }    n\leq N. 
\end{eqnarray}
The function $h^{R_k}$ from \eqref{2.27} converges to zero weakly in $L^2(\Omega)$ because  
$\mathop{\rm supp}\nolimits h^{R_k}\subset \overline{\Omega(R_k)}\setminus \Omega(R_k-1)$ and $\|h^{R_k};L^2(\Omega)\|\leq c(\Omega)$.
The function $g^{R_k}$ is supported in $\overline{\Omega(1)}\setminus\Omega(0)$ and, in view of \eqref{conv3}, uniformly converges to 
$[\Delta,\chi_1]w^\infty$ where 
$$w^\infty(x)=a_n^\infty(e^{-\alpha_n^\infty z_n}+c_n^\infty e^{\alpha_n^\infty z_n})\Phi_1^n(y_n)  \mbox{ in }   Q_n.$$
Note that there appear three options:

1) $a_n^\infty=0 \, \Longrightarrow \, w^\infty(x)=0$ for $x\in Q_n$;

2) $a_n^\infty\ne 0$ and $\alpha_n^\infty\ne 0  \, \Longrightarrow \, c_n^\infty=0$ and  
$w^\infty(x)=a_n^\infty e^{-\alpha_n^\infty z_n}\Phi_1^n(y_n)$ for $x\in Q_n$;

3) $a_n^\infty\ne 0$ and  $\alpha_n^\infty=0  \, \Longrightarrow \,  w^\infty(x)= a_n^\infty(1+c_n^\infty)\Phi_1^n(y_n)$ for $x\in Q_n$.

The function $\widehat{v}^{\,\infty}$ is a solution of the problem
$$
( \nabla\widehat{v}^{\,\infty}, \nabla \psi)_\Omega - \mu^\infty (\widehat{v}^{\,\infty}, \psi)_\Omega= ([\Delta,\chi_1]w^\infty,\psi)_{\Omega} \quad 
\forall\psi\in C^\infty_0(\Omega)
$$
and, therefore, $v^\infty=\widehat{v}^{\,\infty}+\chi_1w^\infty$ becomes a bounded solution of 
the problem \eqref{1.1} with the $\lambda=\mu^\infty$.
Taking into account formula \eqref{2.25} together with relation \eqref{conv2} and using that $\mathcal{X}^R=1$ on $\Omega(R)$ for $R>2$, we obtain  
\begin{multline*}
\|v^\infty;L^2(\Omega(2))\|=\|\widehat{v}^{\,\infty} +\chi_1w^\infty;L^2(\Omega(2))\|=
\\=
\lim_{R_k\to+\infty} \|\widehat{v}^{R_k}+\chi_1w^{R_k};L^2(\Omega(2))\|
=
\lim_{R_k\to+\infty}\|v^{R_k};L^2(\Omega(2))\|=1\,.
\end{multline*}
Thus, $v^\infty\ne 0$.

If $\mu^\infty<\lambda_\dagger$, then $\mu^\infty$ becomes $(\varkappa+1)$-th eigenvalue of the problem \eqref{1.1} that 
contradicts our assumptions.
If $\mu^\infty=\lambda_\dagger$ we obtain the desired result.

Now we are in position to formulate the obtained assertion.

\begin{proposition}
If $\mu^R=\mu_{\varkappa+1}^R < \lambda_\dagger$ for all $R\geq 0$, then there exists threshold resonance in the problem \eqref{1.1}.
\end{proposition}

Propositions 2.2 and 2.6 readily lead to Theorem \ref{Crit1}.

\subsection{A sufficient condition for the absence of trapped modes}
Let us assume that 
\begin{equation}
\label{bfn}
\lambda_\dagger=\Lambda_1^1=\Lambda_1^2=\ldots=\Lambda_1^{\bf n}<\Lambda_1^{\bf n+1}\leq\ldots\leq \Lambda_1^N.
\end{equation}
A characteristic feature of a trapped mode $u\in H^1_0(\Omega)$ looks as follows: 
\begin{equation}
\label{orth}
\int_{\omega_n}u(y_n,0)\Phi_1^n(y_n)dy_n=0,\quad n=1,\ldots, {\bf n}.
\end{equation}
These equalities are supported by the orthogonality conditions in \eqref{A4} and the absence of 
the term $C_n \Phi_1^n(y_n)$ in the Fourier series of the decaying solution $u$ in the outlet $Q_n$.

Let us consider the spectral problem: to find an eigenpair $\{\mu, v^0\}\in\mathbb{R}\times H^1_0(\Omega(0),\Gamma(0))_\perp$
such that
\begin{equation}
\label{ort2}
(\nabla v^0, \nabla \psi)_{\Omega(0)}=\mu (v^0,\psi) \quad \forall \psi\in H^1_0(\Omega(0),\Gamma(0))_\perp.
\end{equation}
Here, $H^1_0(\Omega(0),\Gamma(0))_\perp$ is a subspace of functions in $H^1(\Omega(0))$ which vanish at the
surface $\Gamma(0)=\partial\Omega(0)\cap\partial\Omega$ and enjoy the orthogonality conditions \eqref{orth}.
The differential formulation of this problem involves the equations \eqref{orth} and 
\begin{eqnarray*}
&&-\Delta v^0(x)=\mu v^0(x), \quad x\in \Omega(0), \quad v^0(x)=0, \quad x\in \Gamma(0),\\
&&\partial_\nu v^0(x)=0, \quad x\in \gamma_n(0), \quad n={\bf n}+1, \ldots, N,\\
&&\partial_\nu v^0(x)=C_n\Phi_1^n(y_n), \quad x\in\gamma_n(0), \quad n=1,\ldots, {\bf n},
\end{eqnarray*}
where the constants $C_1,\ldots, C_{\bf n}$ are unfixed.

\begin{theorem}
Let $u$ be a bounded solution of the problem \eqref{1.1} at the threshold $\lambda=\lambda_\dagger$. If the first
eigenvalue of the problem \eqref{ort2} enjoys the relation $\mu_1>\lambda_\dagger$, then $u$ does not decay
at infinity and, therefore, is nothing but a non-trivial stabilizing solution. 
\end{theorem}

\noindent{\bf Proof.}
By the theorem on unique continuation, 
$u$ cannot vanish everywhere in $\Omega(0)$ and, hence, the Friedrichs inequality serving for the problem \eqref{ort2}
gives us the formula
$$
\|\nabla u; L^2(\Omega(0)))\|^2\geq \mu_1\|u; L^2(\Omega(0))\|^2>\lambda_\dagger\|u;L^2(\Omega(0))\|^2.
$$
Taking \eqref{fri} with $R=0$ into account, we come across a contradiction with the integral identity \eqref{1.2}
where $\psi=u$.
\hfill$\Box$

\subsection{Remarks on some known examples}
The papers \cite{na566} and \cite{na578} deal with the symmetric ${\sf T}$- and ${\sf Y}$-shaped planar quantum waveguides where 
multiplicity of the discrete spectrum is 1 while the second eigenvalue of the problem \eqref{1.4} in the smallest node
$\Omega(0)$, the unit square $\square$ and the equilateral triangle $\triangle$ (shaded in Fig. 3, a and b), 
respectively, is strictly bigger than $\lambda_\dagger=\pi^2$.
In this way, the simplest $(\varkappa=1)$ version of Proposition \ref{Pr2.2} applies.

\begin{figure}[t]
\begin{center}
\includegraphics[scale=0.4]{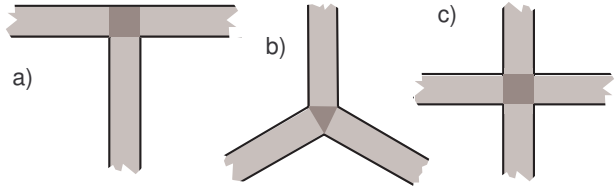}
\end{center}
\caption{{\sf T}-, {\sf Y}- and {\sf X}-shaped waveguides}
\end{figure}

Considering the cruciform waveguide composed from unit circular cylinders, perpendicular to each other, the paper \cite{na602}
demonstrate that $\varkappa=1$ and the eigenvalue $\lambda_2^R$ of the problem \eqref{1.4} with a big $R$ satisfies the inequality
$\lambda_2^R>\lambda_\dagger$, cf. Proposition 2.2. However, for the planar cruciform waveguide made from two perpendicular unit strips, the
Neumann problem in the square $\square$ (shaded in Fig. 3, c) has the eigenvalues $\lambda_1^0=0$, $\lambda_2^0=\lambda_3^0=\pi^2=\lambda_\dagger$. In \cite[\S4]{naIZV}
and \cite[\S 3]{Pank} certain symmetrization tricks were proposed to reject the threshold resonance. At the same time, Proposition 2.5 shows that 
$\lambda_2^R>\pi^2$ when $R>R_2>0$.  
 
\section{An existence criterion}
\subsection{The Steklov-Poincare operator}

To turn the problem \eqref{1.1} with the threshold spectral parameter $\lambda=\lambda_\dagger$ in the infinite domain $\Omega$
into a problem posed in a finite domain, the Steklov--Poincare\footnote{It is also called the Dirichlet-to-Neumann mapping due to its performance.} operator, cf. \cite{Rus, Eng}, is often used. It is expressed 
through solutions of the Dirichlet problem in the semi-infinite cylinder
\begin{eqnarray}
\label{A1}
&&-\Delta U^n=\lambda_\dagger U^n  \mbox{ in }  Q_n,\quad  U^n=0   \mbox{ on }   \Gamma_n:=\partial Q_n\cap \partial\Omega,\\
&& U^n(y,0)=F^n(y)  \mbox{ for }   y\in \omega_n \nonumber.
\end{eqnarray}
Traditionally, this operator acts as follows: $F^n\mapsto \partial_z U^n\big|_{z=0}$.

The Fourier method provides an explicit solution of \eqref{A1} so that the operator takes form
\begin{equation}
\label{def}
\sum_{p=1}^\infty a_p^n \Phi_p^n=F^n\mapsto T^n F^n=-\sum_{p=1}^\infty \kappa_p^n a_p^n\Phi_p^n,
\end{equation}
where $\kappa_p^n=(\Lambda_p^n-\lambda_\dagger)^{1/2}>0$ for $\Lambda_p^n>\lambda_\dagger$
but $\kappa_1^n=-i$ for $n=1,\ldots,{\bf n}$, i.e., in the case $\Lambda_1^n=\lambda_\dagger$ (see \cite[\S2]{na365} for the latter).

If $n={\bf n}+1,\ldots, N$ and $\lambda_\dagger$ stays below the continuous spectrum $[\Lambda_1^n,+\infty)$
of the problem in $Q_n$, and 
\begin{equation}
\label{A6}
(T^nF^n)(y)=\partial_z U^n(y,0),
\end{equation}
where $U^n\in H^1(Q_n)$, is the unique solution of \eqref{A1} with the finite Dirichlet integral. In the case
$n=1,\ldots, {\bf n}$ formula \eqref{A6} is still valid but $U^n$ is a solution to the problem \eqref{A1}
with proper threshold radiation conditions, see Remark \ref{TRC}.

The Fourier method shows that the mapping $T^n: H^1_0(\omega_n)\to L^2(\omega_n)$ is continuous. At the same time,
\begin{equation}
\label{A7}
(T^nF^n,G^n)_{\omega_n}=-\sum_{p=1}^\infty \kappa^n_p a^n_p\overline{b^n_p},
\end{equation}
where $\{a^n_p\}$ and $\{b^n_p\}$ are the Fourier coefficients of $F^n$ and $G^n$, respectively. For $n={\bf n}+1,\ldots, N$,
the relation \eqref{A7} recognizes $T^n$ as a negative operator in the Hilbert space $H^{1/2}_{00}(\omega_n)$,
see \cite[\S1.11]{LiMa}, with the norm 
\begin{equation}
\label{A8}
\|\Psi;H^{1/2}_{00}(\omega_n)\|=(\|\Psi;H^{1/2}(\omega_n)\|^2+\|\rho^{-1/2}\Psi;L^2(\omega_n)\|^2)^{1/2},
\end{equation}
where $\rho=\mathop{\rm dist}\nolimits(y,\partial \omega_n)$ and $H^{1/2}(\omega_n)$ is the standard Sobolev-Slobodetskii space. Notice
that the last weighted norm in \eqref{A8} originates in the Dirichlet condition on $\partial \omega_n$
for the eigenfunctions $\Phi^n_p$. The operator $T^n$ with $n=1,\ldots, {\bf n}$ gets
a skew-symmetric component on the one-dimensional subspace $\mathcal{L}^n$ spanned over the first eigenfunction $\Phi_1^n$
of the problem \eqref{A2}.

Eventually, in the case of the source term $f\in L^2(\Omega)$ with $\mathop{\rm supp}\nolimits f\subset \overline{\Omega(0)}$ a solution of the problem
\begin{eqnarray}
\label{A9}
&& -\Delta u^0-\lambda_\dagger u^0=f   \mbox{ in }  \Omega(0), \quad u^0=0  \mbox{ on }   \Gamma(0),\\
&& \partial_zu^0=T^n u^0   \mbox{ on }   \omega_n(0),\quad n=1,\ldots, N. \nonumber
\end{eqnarray} 
is nothing but the restriction on $\Omega(0)$ of a solution of the problem
\begin{equation}
\label{AAA}
-\Delta u-\lambda_\dagger u=f   \mbox{ in }   \Omega,\quad u=0   \mbox{ on }  \partial \Omega
\end{equation}
with the threshold radiation conditions \eqref{trc}.

\subsection{Symmetrization of the Steklov-Poicare operator}
As was mentioned above, the problem \eqref{A9} inherits all properties of the problem \eqref{AAA}, in particular,
it becomes uniquely solvable if and only if the same property is attributed to \eqref{AAA}.
However, a convenient application of the reduced problem in $\Omega(0)$ needs its unique solvability 
which is clearly absent in the presence of the threshold resonance. In this way, it was proposed in \cite{na365}
to introduce the positive definite symmetric operator
\begin{equation}
\label{A10}
F^n\mapsto M^nF^n =\sum_{p=1}^\infty |\kappa_p^n| a_p^n \Phi_p^n
\end{equation}
and consider the auxiliary problem
\begin{eqnarray}
\label{A11}
&& -\Delta w -\lambda_\dagger w=0   \mbox{ in }   \Omega(0), \quad w=0   \mbox{ on }  \Gamma(0),\\
&& \partial_z w-iM^n w=g^n   \mbox{ on }  \omega_n(0), \quad n=1,\ldots, N\nonumber.
\end{eqnarray}
The weak formulation of this problem reads: to find $w\in \mathcal{H}^0_0$, see Section 2.2,
such that 
\begin{equation}
\label{A12}
(\nabla w,\nabla v)_{\Omega(0)}-\lambda_\dagger (w,v)_{\Omega(0)} - i\langle {\bf M} w,v\rangle=\langle g,v\rangle
\quad \forall v\in \mathcal{H}^0_0.
\end{equation}
Here, ${\bf M}=\mathop{\rm diag}\nolimits\{M^1,M^2,\ldots, M^N\}$, $g=(g^1,g^2,\ldots, g^N)$ and $\langle\,,\rangle$ is the extension
of the scalar product in ${\bf L}:= L^2(\omega_1(0))\oplus \ldots\oplus L^2(\omega_n(0))$ 
up to the duality 
between the space 
$$
{\bf H}=H^{1/2}_{00}(\omega_1(0))\oplus \ldots \oplus H^{1/2}_{00}(\omega_N(0))
$$
and its adjoint ${\bf H}^*:=H^{-1/2}_{00}(\omega_1(0))\oplus \ldots \oplus H^{-1/2}_{00}(\omega_N(0))$.

As was proved in \cite[Lemma 2.2]{na365} and can be easily verified on the basis
of the theorem on unique continuation, in view of the presence of the skew-symmetric sesquilinear 
form $i\langle {\bf M} w,v\rangle$  
the problem \eqref{A12} with any $g\in {\bf H}^*$
has a unique solution $w\in \mathcal{H}^0_0$ and the following estimate is valid:
\begin{equation}
\label{A13}
\|w;\mathcal{H}^0_0\|\leq c\|g;{\bf H}^*\|.
\end{equation}

\subsection{The fictitious scattering operator.} 
Following \cite{na365}, we introduce an artificial object, a unitary operator ${\bf S}$ in ${\bf L}$ which
can be directly constructed through solutions of the uniquely solvable problem \eqref{A11} and becomes
an identificator of all bounded solutions at the threshold, see Theorem \ref{Th34}.
 
Let ${\bf M}^{1/2}$ be the positive square root of the positive self-adjoint operator ${\bf M}$ in \eqref{A10}.
For any $\psi\in {\bf L}$, 
we denote by $w(\psi)\in \mathcal{H}^0_0$ the (unique) solution of the problem 
\eqref{A12} with the specific right-hand side
\begin{equation}
\label{A14}
g=-2^{1/2} i {\bf M}^{1/2} \psi \in {\bf H}^*
\end{equation}
and set 
\begin{equation}
\label{A15}
{\bf S} \psi=i\psi-2^{1/2}i{\bf M}^{1/2} w(\psi)|_{\omega(0)} \in {\bf L}\,,
\end{equation}
where $\omega(0)=\omega_1(0)\times\ldots\times\omega_N(0)$.
In view of the estimate \eqref{A13} and the properties of the operator ${\bf M}$ we see that
\eqref{A15} is a continuous operator in ${\bf L}$. Moreover, in \cite[Theorem 2.1]{na365}
it is verified that, owing to the special choice \eqref{A14} of the right-hand side in \eqref{A11}, 
${\bf S}$ is a unitary operator in ${\bf L}$.

\subsection{The criterion for trapped modes}
Let ${\bf L}_0$ be the subspace
\begin{eqnarray}
\label{A16}
\{\psi\in{\bf L} :&& \psi|_{\omega_n(0)}=c_n \Phi_1^n, n=1,\ldots, {\bf n}, c_n\in\mathbb{C},\\
&&\psi|_{\omega_n(0)}=0, n={\bf n}+1,\ldots,N\}\nonumber
\end{eqnarray}
and let ${\bf L}_\perp ={\bf L} \ominus {\bf L}_0$ be the orthogonal complement of \eqref{A16}.
Denoting the orthogonal projectors on ${\bf L}_0$ and ${\bf L}_\perp$ by ${\bf P}_0$ and ${\bf P}_\perp$,
respectively, we define the operator
\begin{equation}
\label{A17}
{\bf S}_\perp={\bf P}_\perp{\bf S} {\bf P}_\perp:{\bf L}_\perp\to {\bf L}_\perp.
\end{equation}

In \cite[Theorem 3.1]{na365} it is verified that the mapping
$$
{\bf D}_{tr}\ni u\mapsto \psi = 2^{-1/2}(1+i){\bf M}^{1/2} u|_{\omega(0)} \in \ker ({\bf S}_\perp-\mathop{\rm Id}\nolimits_\perp)
$$
is a bijection where ${\bf D}_{tr}$ is the subspace of trapped modes in the problem \eqref{1.1}
at the threshold $\lambda=\lambda_\dagger$ and $\ker ({\bf S}_\perp-\mathop{\rm Id}\nolimits_\perp)$
is the eigenspace of the operator \eqref{A17} for its eigenvalue 1. This fact readily
establishes the existence criterion for trapped modes.

\begin{proposition} 
There holds
\begin{equation}
\label{tr}
\dim {\bf D}_{tr}=\dim \ker ({\bf S}_\perp-\mathop{\rm Id}\nolimits_\perp),
\end{equation}
i.e. a trapped mode exists if and only if the operator \eqref{A17} has the eigenvalue~1.
\end{proposition}

It should be mentioned that 
\begin{equation}
\label{A19}
\begin{array}{l}
\psi\in \ker ({\bf S}_\perp-\mathop{\rm Id}\nolimits_\perp)\quad \Rightarrow \quad \|\psi;{\bf L}\|^2=\|{\bf S}\psi;{\bf L}\|^2=\|{\bf P}_\perp {\bf S}\psi;{\bf L}\|^2+\\
+\|{\bf P}_0{\bf S} \psi;{\bf L}\|^2
=\|\psi;{\bf L}\|^2+\|{\bf P}{\bf S} \psi;{\bf L}\|^2\quad \Rightarrow \quad {\bf P}{\bf S}\psi=0.
\end{array}
\end{equation} 
In other words, $\psi\in \ker({\bf S}_\perp -\mathop{\rm Id}\nolimits_\perp)$ is an eigenfunction of the intact fictitious scattering operator ${\bf S}$
corresponding to the eigenvalue 1.
 
\subsection{Threshold radiation conditions and the threshold scattering matrix}
At the threshold $\lambda_\dagger$ the standing $\Phi_1^n(y_n)$ and resonance $y_n\Phi_1^n(y_n)$ waves occur in the outlets 
$Q_n$, $n=1,\ldots, {\bf n}$. These waves cannot be classified by classical Sommerfeld radiation principle because of their null wave
number. In order to define a unitary and symmetric scattering matrix at the threshold, we follow \cite[Ch.5, \S3]{NaPl},
and introduce the couples of linear in $z_n$ waves
\begin{equation}
\label{Z1}
w_n^{in}(x)=\chi(z_n)2^{-1/2}(z_n+i)\Phi^n_1(y_n), \quad w_n^{out}(x)=\chi(z_n)2^{-1/2}(z_n-i)\Phi_1^n(y_n)
\end{equation} 
where the superscripts mean ``incoming'' and ``outgoing''. The linear combinations \eqref{Z1}
of the resonance and standing waves emerging at the threshold possess the remarkable properties:
\begin{equation}
\label{Z2}
w_n^{in}(x)=\overline{w_n^{out}(x)}
\end{equation} 
and
\begin{equation}
\label{Z3}
\begin{array}{l}
q_n(w_n^{in},w_n^{in})=-i, \quad q_n(w_n^{out},w_n^{out})=i,\\
q_j(w_n^{in},w_n^{out})=-\overline{q_n(w_n^{out},w_n^{in})}=0
\end{array}
\end{equation}
with the sesquilinear and anti-Hermitian form
\begin{equation}
\label{Z4}
q_n(u,v)=\int_{\omega_n}\big( \overline{v(x)}\partial_{z_n}u(x)-u(x)\overline{\partial_{z_n}v(x)}\big)\big|_{z_n=R}dy_n
\end{equation} 
which appears as a surface integral in the Green formula on the truncated waveguide \eqref{1.0} and, therefore,
is independent of $R>1$ for waves \eqref{Z1} and their linear combinations.

\begin{remark}
\label{TRC}
The threshold radiation condition for the problem \eqref{A9} reads
\begin{equation}
\label{trc}
u-\sum_{n=1}^{\bf n} c_n w_n^{out}\in H^1(\Omega)
\end{equation}
where ${\bf n}$ is defined in \eqref{bfn}, $w_n^{out}$ is the outgoing wave in \eqref{Z1}
and $c_1,\ldots, c_{\bf n}$ are some coefficients. Conditions of type \eqref{trc}
have been introduced in \cite[Ch. 5]{NaPl}, as well as their straight-forward
modifications for the threshold inside the continuous spectrum (the eigenvalues $\Lambda_p^n$ of the model
problem \eqref{A2} with $p\geq 2$). The corresponding problems always inherit all important properties of the problems 
outside the thresholds.
\end{remark}

As was demonstrated in \cite[\S3 Ch.5]{NaPl}  and, e.g., \cite{na489}, the relation \eqref{Z3} and \eqref{Z2}
are sufficient to guarantee the existence of the special solutions 
\begin{equation}
\label{Z5}
Z_n(x)=w_n^{in}(x)+\sum_{k=1}^{{\bf n}} {\bf s}_{kn} w_k^{out}(x)+\widetilde{Z}_n(x)
\end{equation}
to the problem \eqref{1.1} with $\lambda=\lambda_\dagger$ as well as the unitary and symmetry properties of the 
threshold scattering matrix ${\bf s}$ composed of the coefficients ${\bf s}_{kn}$, $k, n=1,\ldots, {\bf n}$, in \eqref{Z5}.
Note that $Z_n(x)$ decays in the outlets $Q_{{\bf n}+1}, \ldots, Q_{N}$ only but the reminder $\widetilde{Z}_n$ does in all outlets.

\begin{remark}
The form \eqref{Z4} induces an indefinite metrics in the $2{\bf n}$-dimensional subspace $\mathcal{W}$ of polynomial waves, and, of course, the 
above-mentioned basis in $\mathcal{W}$ is not unique. For example, the waves 
\begin{equation}
\label{Z6}
\begin{array}{l}
{\bf w}_n^{in}(x)=\chi(z_n)2^{-1/2}e^{i\psi}(1-iz_n)\Phi_1^n(y_n),\\
{\bf w}_n^{out}(x)=\chi(z_n)2^{-1/2}e^{-i\psi}(1+iz_n)\Phi_1^n(y_n)
\end{array}
\end{equation}
with $\psi\in\mathbb{R}$ verify the same relations \eqref{Z2} and \eqref{Z3} as waves \eqref{Z1}. The threshold scattering matrix ${\bf s}$ initiated by incoming
waves in \eqref{Z6} is equal to $e^{2i\psi}{\bf s}$. This observation will allow us to formulate in Theorem 3.4 the
common criterion for the existence of trapped modes and stabilizing solutions.
\end{remark} 
 
\subsection{The criterion for the existence of stabilizing solutions.}
The following assertion can be found in the paper \cite{na566} but its proof is very simple and we 
reproduce it here for reader's convenience. We also mention that other arguments in \cite{MoVa} and \cite{Gri}
had let to similar assertions expressed in different terms.

\begin{proposition}
\label{Pr1}
Dimension of the subspace ${\bf D}_{st}$ of stabilizing solutions coincides with multiplicity of the eigenvalue $-1$ of 
the threshold scattering matrix ${\bf s}$, i.e. $\dim {\bf D}_{st}=\dim\ker ({\bf s}+\mathop{\rm Id}\nolimits_0)$, where $\mathop{\rm Id}\nolimits_0$ is the unit matrix of size ${\bf n}\times{\bf n}$. 
If ${\bf s} c+c=0$ for a column $c\in \mathbb{C}^{\bf n}\setminus\{0\}$, then a nontrivial stabilizing solution is given by the linear combination
\begin{equation}
\label{Z7}
Z=c_1Z_1+\ldots+c_{\bf n} Z_{\bf n}\,.
\end{equation}
\end{proposition}

\noindent{\bf Proof.}
The function \eqref{Z7} admits the decomposition
\begin{equation}
{\displaystyle
\label{Z8}
\begin{array}{l}
Z(x)=\sum\limits_{n=1}^{{\bf n}}\Big(c_n w_n^{in}(x)+\sum\limits_{k=1}^{\bf n} c_k {\bf s}_{kn} w_n^{out}(x)\Big)+\widetilde{Z}(x)=\\
\sum\limits_{n=1}^{\bf n} c_n(w_n^{in}(x)-w_n^{out}(x))+\widetilde{Z}(x)=-\sqrt{2}i\sum\limits_{n=1}^{\bf n} c_n \chi_n(z_n)\Phi_1^n(y_n)
+\widetilde{Z}(x)\,.
\end{array}
}
\end{equation}
Here, we used the equality ${\bf s} c+c=0\in \mathbb{C}^{\bf n}$ and formulas \eqref{Z1} to observe that $Z$ is bounded and does not decay at infinity.
Reading the chain \eqref{Z8} from right to left proves the equalities $c_n=-\sum c_k {\bf s}_{kn}$, $n=1,\ldots, {\bf n}$, and concludes
with the whole assertion. 
\hfill$\Box$ 

In other words, the threshold scattering matrix contain the complete information on stabilizing solutions of the problem \eqref{1.1}
with $\lambda=\lambda_\dagger$.

\subsection{The fictitious scattering operator and stabilizing solutions}
The function $-iZ_n$, see \eqref{Z5}, satisfies the problem \eqref{A11} with the right-hand side
$$
g^n=-i\partial_z Z_n-M^nZ_n   \mbox{ on }   \omega_n, \quad n=1,\ldots, {\bf n}.
$$
Since ${\bf M}^{1/2}{\bf P}_0={\bf P}_0{\bf M}^{1/2}={\bf P}_0$ according to definitions \eqref{A10} and \eqref{def}, we take \eqref{Z5} into account and obtain
\begin{equation}
\label{X1}
{\bf P}_0g^n=-2^{1/2}i{\bf e}_n \Phi_1^n \quad \mbox{where} \quad {\bf e}_n=(\delta_{1,n},\delta_{2,n},\ldots, \delta_{N,n}).
\end{equation}
Comparing \eqref{X1} with \eqref{A14} and recalling \eqref{A16} yield
$$
{\bf P}_0{\bf S}{\bf P}_0\psi=i\psi-2^{1/2}{\bf M}^{1/2}{\bf P}_0 Z|_{\omega(0)}=i\psi-2^{1/2}(i2^{-1/2}\psi-i2^{-1/2}{\bf s}\psi)=i{\bf s}\psi
$$
where $Z$ is the linear combination \eqref{Z7} and
$$
\psi=(c_1\Phi_1^1,\ldots, c_{\bf n} \Phi_1^{\bf n},0,\ldots,0),\quad {\bf s}\psi=\sum_{k=1}^nc_k ({\bf s}_{1k}\Phi_1^1,\ldots, {\bf s}_{{\bf n} k} \Phi_1^{\bf n},0,\ldots,0).
$$
In other words, the operator
\begin{equation}
\label{X2}
{\bf S}_0={\bf P}_0{\bf S}{\bf P}_0:{\bf L}_0\to{\bf L}_0
\end{equation}
realizes as the unitary matrix $i{\bf s}$ that allows us to reformulate the criterion in Proposition 3.3 in terms of the operator \eqref{X2}, 
namely
\begin{equation}
\label{X3}
\dim {\bf D}_{st}=\dim\ker ({\bf S}_0+i\mathop{\rm Id}\nolimits_0).
\end{equation}

Repeating the calculations \eqref{A19} we see that ${\bf P}_\perp{\bf S}\psi=0$ in the case $\psi\in\ker ({\bf S}_0+i\mathop{\rm Id}\nolimits_0)$ and, therefore, 
$\psi\in \ker ({\bf S}+i\mathop{\rm Id}\nolimits)$. Thus,
formulas \eqref{tr} and \eqref{X3} lead to the following criterion for the existence of bounded solutions of the problem \eqref{1.1}
with $\lambda=\lambda_\dagger$, that is, for the threshold resonance.

\begin{theorem}
\label{Th34}
The subspace ${\bf D}_{bd}={\bf D}_{tr}\cup {\bf D}_{st}$ of bounded solutions verifies the relation
$$
\dim{\bf D}_{bd}=\dim\ker (\widehat{{\bf S}}-\mathop{\rm Id}\nolimits),
$$
where $\mathop{\rm Id}\nolimits$ is the identify operator in ${\bf L}$ and
\begin{equation}
\label{X4}
\widehat{{\bf S}}=({\bf P}_\perp+2^{-1/2}(1-i){\bf P}_0){\bf S} ({\bf P}_\perp +2^{-1/2}(1-i){\bf P}_0).
\end{equation}
\end{theorem}

We emphasize that operator \eqref{X4} is still unitary.

\section{Acknowledgments}
Research is financially supported by grant № 17-11-01003 of the Russian Science Foundation.

\end{document}